\documentclass[11pt,twoside]{article}
\renewcommand\baselinestretch{1.3}%
\newcommand\kaishu{\CJKfamily{kai}}
\newcommand\heiti{\CJKfamily{hei}}
\newcommand\songti{\CJKfamily{song}}

\newcommand{\sihao}{\fontsize{14pt}{\baselineskip}\selectfont}

\newcommand{\wuhao}{\fontsize{10.5pt}{\baselineskip}\selectfont}
\newcommand{\xiaowuhao}{\fontsize{9pt}{\baselineskip}\selectfont}

\usepackage{CJK,indentfirst,amsmath,amsfonts,amssymb,amsthm,cite}
\usepackage[dvips]{graphics} %
\usepackage{graphicx}
\usepackage[dvips]{color} %
\usepackage{cite}
\usepackage{fancyhdr} %
\usepackage{ifthen} %
\usepackage{cases}
\begin{document}
\begin{CJK*}{GBK}{song}
\renewcommand\baselinestretch{1.2}
\renewcommand\arraystretch{1.2}
\catcode`@=11 \long\def\@makefntext#1{\parindent 2em\indent \hbox
to 0pt{\hss$^{}$}#1} \catcode`\@=12
\pagestyle{fancy}%
\setcounter{page}{1} %
\newcounter{jie}
\renewenvironment{proof}[1][\heiti 证明]{\textbf{#1.} }{
\begin{flushright}
$\Box$
\end{flushright}
}
%
\newboolean{first}%
\setboolean{first}{true}%
\pagestyle{fancy}
\fancyhf{}   %
\fancyhead[RE]{\xiaowuhao{\songti {}\hspace{2pt}}}
\fancyhead[LO]{\xiaowuhao{\ \songti {}}} \fancyhead[LE]{\ \thepage}
\fancyhead[RO]{\thepage \hspace{2pt}}
\fancyhead[CE]{\xiaowuhao{\songti   }}
\fancyhead[CO]{\xiaowuhao{\songti   On lower bounds of the first
eigenvalue of Finsler-Laplacian}}
%

\newcommand{\makefirstpageheadrule}{%
\makebox[0pt][l]{\rule[0.55\baselineskip]{\headwidth}{0.4pt}}%
\rule[0.7\baselineskip]{\headwidth}{0.4pt}}

\newcommand{\makeheadrule}{%
\rule[0.7\baselineskip]{\headwidth}{0.4pt}}
\renewcommand{\headrule}{%
\ifthenelse{\boolean{first}}{\makeheadrule}{\makefirstpageheadrule}%
}%

\title{    On  Lower Bounds of the First Eigenvalue  of Finsler-Laplacian }
\vspace{0.6cm}
\author{\wuhao\songti  Songting Yin,  Qun He and Yi-Bing Shen \\
}   %

\footnotetext[1]{
\begin{minipage}[t]{\textwidth}
{
}
\end{minipage}
}

\date{}
\maketitle
\renewcommand\baselinestretch{1.2}

\vspace*{-0.5cm} {
\begin{minipage}{\textwidth}
\begin{minipage}[t]{14cm}
{\wuhao\heiti  \textbf{Abstract} \wuhao\kaishu \ \ By using Bochner
technique and gradient estimate, we give the lower bound estimates of
the first eigenvalue of Finsler-Laplacian on Finsler manifolds.
These results generalize the corresponding famous theorems in the
Riemannian geometry.
 }

\end{minipage}

\ \

\end{minipage}
} \vspace*{0.2cm}
\section*{\sihao\heiti 1.  \ Introduction} %
The research on the first (nonzero) eigenvalue of Laplacian  plays
an important role in global differential geometry. In the Riemannian
case, Lichnerowicz [10] advocated it for the first time and gave the
lower bound estimate of the first eigenvalue via the restriction of
the Ricci curvature. Afterwards, Obata [12] further  established a
rigidity theorem, demonstrating the optimality of Lichnerowicz'
estimate. For the non-negative Ricci curvature, Li-Yau [9] employed
the gradient estimates of the eigenfunctions and got the lower bound
estimate of the first eigenvalue via the diameter of the manifolds.
Then this method was improved further and the optimal result was
obtained by Zhong-Yang [22]. Recently, Hang-Wang [7] proved that $S^{1}$
is the only case for the first eigenvalue attaining its lower bound.
 Precisely, they  achieved the following results respectively.\\

 \textbf{\wuhao\heiti Theorem 1.1. ([Li]-[Ob])} \ \  \ \emph{Let $(M,g)$ be an $n-$dimensional compact Riemannian manifold
without boundary. If the Ricci curvature satisfies}
$$\textmd{Ric}_{M}\geq (n-1)k$$
\emph{for some given constant $k>0$, then}
$$\lambda_{1}\geq nk,$$
\emph{where the equality holds if and only if $M$ is isometric to the n-sphere of constant sectional curvature k, so that
the diameter of M is } $\frac{\pi}{\sqrt{k}}.$\\

\textbf{\wuhao\heiti Theorem 1.2.([LY]-[ZY]-[HW])} \ \  \ \emph{Let
$(M,g)$ be an $n-$dimensional compact Riemannian manifold without
boundary. If }$\textmd{Ric}_{M}\geq 0$, \emph{then
$$\lambda_{1}\geq \frac{\pi^{2}}{d^{2}},$$
where d denotes the diameter of $(M,g)$ and the
 equality holds if and only if $M$ is isometric to $S^{1}(\frac{d}{\pi})$}.\\

As a natural generalization of Riemannian manifolds, Finsler
manifolds are differentiable manifolds of which on each tangent
space one endows a Minkowski norm instead of a Euclidean norm.
Recent studies on Finsler manifolds have taken on a new look.  Up to
now, there have been several different definitions of
Finsler-Laplacians, introduced respectively by Bao-Lacky [3],
Antonelli-Zastawniak [1], Centro [4], Thomas [19], and Ge-Shen [6].
By using the Finsler-Laplacian, Y. Ge and Z. Shen  gave the
Faber-Krahn type inequality for the first Dirichlet eigenvalue of
the Finsler-Laplacian in [6]. B. Wu and Y. Xin [20] proved that for
a complete noncompact and simply connected Finsler manifold with
finite reversibility $\lambda$ and nonpositive flag curvature, if
$\textmd{Ric}\leq-a^{2} (a>0)$ and $\sup_{M}\|S\|<a$, then
$\lambda_{1}\geq\frac{(a-\sup_{M}\|S\|)^{2}}{4\lambda^{2}}$. Another
interesting result on this direction, due to
 G. Wang and C. Xia [21], says  that for a compact Finsler measure space,
 if the \emph{weighted Ricci curvature} (see Definition 2.1 below) $\textmd{Ric}_{N}\geq K$, $N\in[n, \infty], K\in R$,
 then $\lambda_{1}\geq\lambda_{1}(K,N,d)$ where $\lambda_{1}(K,N,d)$ represents the first eigenvalue of
 the 1-dimensional problem (see [21] for details).

In this paper we focus on lower bound estimates of the  first
eigenvalue of the Finsler-Laplacian$^{[6]}$ on Finsler manifolds
with an arbitrary volume form $d\mu $.  The main purpose is to
generalize Theorem 1.1 and Theorem 1.2 into the Finsler case. It
should be noted that since the Finsler-Laplacian is a nonlinear
operator, some methods  used in the Riemannian case are not
adaptable any more. To overcome these difficulties, we have to
utilize the properties of the {weighted gradient} and the {weighted
Laplacian} in weighted Riemannian manifold $(M,g_{V})$$^{[13, 21]}$.
Here the weighted gradient and weighted Laplacian play an important
and reasonable role in studying the first eigenvalue of the
Finsler-Laplacian. With the help of them,  we can convert some
nonlinear problems into the linear ones and even calculate something
as simple as in the Riemannian case. In addition, we also make use
of Bochner technique and some gradient estimates. Then by using the
restriction of weighted-Ricci curvature$^{[14]}$  (in the Riemannian
case it is just the Ricci curvature), we obtain the estimates on the
lower bounds
for the first eigenvalue $\lambda_1$ of Finsler-Laplacian in $(M,F,d\mu )$. Concretely, we get the main results
as follows. \\

\textbf{\wuhao\heiti Theorem 3.1.} \ \ \ \emph{Let $(M,F)$ be an $n-$dimensional forward geodesically complete
connected Finsler manifold. If the weighted Ricci curvature and $S$-curvature satisfy}
$$\textmd{Ric}_{N}\geq (n-1)k,~~~~~~~~~~~~\dot{S}\leq \frac{(N-n)(n-1)}{N-1}k$$
\emph{for some uniform positive constant $k$ and $N\in (n,\infty)$,
where $\dot{S}$ denotes the change rate of the $S$-curvature along
geodesics, then
$$\lambda_{1}\geq \frac{n-1}{N-1}Nk.$$
Moreover, the diameter of $M$ is $\sqrt{\frac{N-1}{n-1}}\frac{\pi}{\sqrt{k}}$ if the equality holds.}\\

\textbf{\wuhao\heiti Theorem 3.3.} \ \ \ \emph{Let $(M,F)$ be an $n-$dimensional forward geodesically complete connected Finsler manifold. If $S=0$ and Ricci curvature  }
$\textmd{Ric} \geq (n-1)k$
\emph{for some uniform positive constant $k$, then
$$\lambda_{1}\geq nk.$$
Moreover, if the equality holds, then the diameter of $M$ is $\frac{\pi}{\sqrt{k}}$,  and $M$ is homeomorphic to $S^{n}$. In particular, if $F$
is reversible and $M$ has Busemann-Hausdorff volume form, then $(M,F)$ is isometric to $S^{n}(\frac{1}{\sqrt{k}})$. }\\

\textbf{\wuhao\heiti Theorem 3.4.} \  \emph{Let $(M,F)$ be an $n-$dimensional compact Finsler manifold. If the weighted Ricci curvature } $\textmd{Ric}_{\infty}\geq 0$,
\emph{then
$$\lambda_{1}\geq\frac{\pi^{2}}{d^{2}},$$
where $d$ denotes the diameter of $(M,F)$.}\\

\ \ Here the term "weighted Ricci curvature"(Definition 2.1) and the
notation "$\dot{S}$" (Definition 2.2) will be given in section 2
below. If $F$ is Riemannian metric,   above Theorems are in accord
with Theorem 1.1 and Theorem 1.2.\\

\vspace*{-0.4cm}
\section*{\sihao\heiti 2.  \ Preliminaries} %
   Throughout this paper, we assume that $M$ is an $n-$dimensional oriented smooth manifold without boundary. A $Finsler$ $metric$ on $M$ is a function $F: TM\longrightarrow[0,\infty)$ satisfying the following properties (i) $F$ is smooth on $TM\backslash{0}$; (ii) $F(x,\lambda y)=\lambda F(x,y)$ for all $\lambda>0$; (iii) the induced quadratic form $g$ is positive-definite, where
         $$g:=g_{ij}dx^{i} \otimes dx^{j}, ~~~~~~~~~   g_{ij}=\frac{1}{2}[F^{2}] _{y^{i}y^{j}}.  $$
   Here and from now on, we will use the following convention of index ranges unless other stated:
              $$1\leq i, j\cdots \leq n ;~~~~~~1\leq \alpha, \beta \cdots \leq n-1,~~~~~~~~~\bar{\alpha}=n+\alpha.$$

   The projection $\pi : TM\longrightarrow M $ gives rise to the pull-back bundle $\pi^{\ast}TM$ and its dual bundle $\pi^{\ast}T^{\ast}M$ over $TM\backslash{0}$. In $\pi^{\ast}T^{\ast}M$ there is a global section $\omega=[F]_{y^{i}}dx^{i},$ called the $Hilbert$ $form$, whose dual is $\ell=\ell^{i}\frac{\partial}{\partial x^{i}},\ell^{i}=\frac{y^{i}}{F}$, called the $distinguished$ $field.$

  Let $\{e_{i}\}_{i=1}^{n}$ be a local orthonormal basis on $\pi^{\ast}TM$ such that its dual basis is $\{\omega^{i}\}_{i=1}^{n}$ with $\omega^{n}=\omega$. As is well known that on the pull-back bundle $\pi^{\ast}TM$ there exists uniquely the Chern connection $^{c}\nabla$ with $^{c}\nabla e_{i}=\omega_{i}^{j}e_{j}$ satisfying
  $$d\omega^{i}=-\omega^{i}_{j}\wedge \omega^{j},~~~~~~~~~~~\omega^{\alpha}_{n}=\omega^{\bar{\alpha}},~~~~~\omega^{n}_{n}=0,$$
  $$\omega^{i}_{j}+\omega^{j}_{i}=-2C_{ijk}\omega^{k}_{n},~~~~~~~~~~~C_{njk}=0,~~~~~~~~~~~~~~~~~~~~~$$
  where $C_{ijk}=\frac{1}{F}A_{ijk}$ is called the Cartan tensor.

Let $u:M\longrightarrow R$ be a smooth function. Then we can view $u$ as its lift on the projective sphere bundle $SM$. Define
$$du:=u_{i}\omega^{i},~~\eqno(2.1)$$
$$du_{i}-u_{j}\omega^{j}_{i}:=u_{i|j}\omega^{j}+u_{i;\alpha}\omega^{\bar{\alpha}},\eqno(2.2)$$
where $"|"$ and $";"$ denote the horizontal covariant derivative with respect to $^{c}\nabla$ and vertical derivative, respectively. Taking exterior differentiation of (2.1) and making use of (2.2), the structure equations with respect to the Chern connection, we have
$$u_{i|j}=u_{j|i},~~~~~~~~~~~~~~~~u_{i;\alpha}=0.$$

  The curvature 2-forms of the Chern connection $^{c}\nabla$ are
  $$d\omega^{i}_{j}-\omega^{k}_{j}\wedge \omega^{i}_{k}=\Omega^{i}_{j}:=\frac{1}{2}R^{i}_{j~kl}\omega^{k}\wedge \omega^{l}+P^{i}_{j~k\alpha}\omega^{k}\wedge \omega^{\bar{\alpha}},$$
  where $R^{i}_{j~kl}=-R^{i}_{j~lk}$ and $P^{i}_{j~k\alpha}=P^{i}_{k~j\alpha}$.
  The Landsberg curvature is defined as $P^{i}_{~jk}:=P^{i}_{n~jk}$, which satisfies
  $$P_{ijk}=\delta_{il}P^{l}_{~jk}=-\dot{A}_{ijk},~~~~~~~~P_{n\alpha\beta}=0,$$
  where "." denotes the covariant derivative along the Hilbert form. The flag curvature tensor can be defined by
  $$R_{\alpha\beta}=\delta_{\alpha\gamma}R^{\gamma}_{n~\beta n}.$$
  For a unit vector $V=V^{i}e_{i}$, the flag curvature $K(y;V)$ is
  $$K(y;V)=R_{\alpha\beta}V^{\alpha}V^{\beta}.$$
  The Ricci curvature for $(M,F)$ is defined as
  $$\textmd{Ric}\triangleq \textmd{Ric}(y)=\sum_{\alpha=1}^{n-1} K(y;e_{\alpha})=\sum_{\alpha=1}^{n-1} R_{\alpha\alpha}.$$

Clearly, Ricci curvature $\textmd{Ric}(y)$ is positively homogeneous of degree zero. i.e., $\textmd{Ric}(\lambda y)=\textmd{Ric}(y)$ for all $\lambda>0$. Now we can introduce the weighted Ricci curvature on the Finsler manifolds, which was defined by Ohta in [14], motivated by the work of Lott-Villani [11] and Sturm [18] on metric measure space.\\

\textbf{\wuhao\heiti Definition 2.1.([14])} \ \  ~~~Given a  vector $V\in T_{x}M$, let $\eta:(-\varepsilon,\varepsilon)\longrightarrow M$ be the geodesic such that $\eta^{'}(0)=V$. We set $d\mu =e^{-\Psi}\textmd{vol}_{\eta^{'}}$ along $\eta$, where $\textmd{vol}_{\eta^{'}}$ is the volume form of $g_{\eta^{'}}$. Define \emph{weighted  Ricci curvature} by

\begin{itemize}
   \item $\textmd{Ric}_{N}(V):=\textmd{Ric}(V)+\frac{(\Psi\circ \eta)^{''}(0)}{F(V)^{2}}-\frac{(\Psi\circ \eta)^{'}(0)^{2}}{(N-n)F(V)^{2}}~~~~~~$\texttt{for}~~ $N\in (n,\infty),$
   \item $\textmd{Ric}_{\infty}(V):=\textmd{Ric}(V)+\frac{(\Psi\circ \eta)^{''}(0)}{F(V)^{2}}.$
   \end{itemize}

\textbf{\wuhao\heiti Remark:}\ \  The above definition is slightly different from that in [14] where the weighted Ricci curvature is  positively homogeneous of degree two.\\

As is well known that $S$-curvature is one of the most important non-Riemannian quantities in Finsler geometry. For any $ y\in T_{x}M\backslash{0}$, let $\gamma(t)$ be the geodesic with $\gamma(0)=x,\dot{\gamma}(0)=y$. Then \emph{$S$-curvature} is defined by
 $$S(x,y)=\frac{d}{dt}[\tau(\gamma(t),\dot{\gamma}(t))]_{t=0}.$$
 An $n$-dimensional Finsler metric $F$ on a manifold is said to have \emph{constant} $S$-curvature if $S=(n+1)cF$ for some constant $c$. In order to  measure the rate of change of the $S$-curvature along geodesics, we give the following\\

\textbf{\wuhao\heiti Definition 2.2.} \ \ \ For any $y\in T_{x}M\backslash{0}$, define
$$\dot{S}(x,y)=\frac{1}{F^{2}}\frac{d}{dt}[S(\gamma(t),\dot{\gamma}(t))]_{t=0},\eqno(2.3)$$
where $\gamma(t)$ is geodesic satisfying $\gamma(0)=x,\dot{\gamma}(0)=y$.\\

 \textbf{\wuhao\heiti Remark:}\ \
By Definition 2.2, we get $ \dot{S}(x, y)=\frac{1}{F^{2}}S_{|i}y^{i}=\frac{1}{F^{2}}\{S_{x^{i}}y^{i}-2S_{y^{i}}G^{i}\}.$
It follows that
 $\dot{S}(x,\lambda y)= \dot{S}(x,y),~\forall\lambda>0.$
In addition, according to Definition 2.1, $d\mu =e^{-\Psi}\textmd{vol}_{\eta^{'}}$ implies $\Psi=\tau$ along geodesic $\eta$, here $\tau$ denotes the \emph{distortion}  of $F$ with respect to $d\mu$. So by definition of $S$ and $\dot{S}$ we have
$$S=(\Psi\circ \eta)^{'}(0),~~~~~~~~~~~ \dot{S}=\frac{(\Psi\circ \eta)^{''}(0)}{F^{2}},\eqno(2.4)$$
where $\Psi, \eta$ are defined in Definition 2.1.\\

Let $X=X^{i}\frac{\partial}{\partial x^{i}}$ be a vector field. Then the \emph{covariant derivative} of $X$  by $v\in T_{x}M$ with reference  vector $w\in T_{x}M\backslash 0$ is defined by
$$D^{w}_{v}X(x):=\left\{v^{j}\frac{\partial X^{i}}{\partial x^{j}}(x)+\Gamma^{i}_{jk}(w)v^{j}X^{k}(x)\right\}\frac{\partial}{\partial x^{i}},\eqno(2.5)$$
where $\Gamma^{i}_{jk}$ denotes the coefficients of the Chern connection given by
$$\Gamma^{i}_{jk}=\frac{1}{2}g^{il}(\frac{\delta g_{lj}}{\delta x^{k}}+\frac{\delta g_{lk}}{\delta x^{j}}-\frac{\delta g_{jk}}{\delta x^{l}}).$$
And $$\frac{\delta}{\delta x^{i}}=\frac{\partial}{\partial x^{i}}-N^{j}_{i}\frac{\partial}{\partial y^{j}},~~~~~~~N^{j}_{i}=\frac{\partial G^{j}}{\partial y^{i}},~~~~~G^{i}=\frac{1}{4}g^{il}\{[F^{2}]_{x^{k}y^{l}}y^{k}-[F^{2}]_{x^{l}}\}.$$\\

Now let $L^{\ast}:T^{\ast}M\longrightarrow TM$ denote the\emph{ Legendre transform}. Then $L^{\ast}$ is norm-preserving map satisfying $L^{\ast}(a\zeta)=aL^{\ast}(\zeta)$, for all $a>0, \zeta\in T^{\ast}M$.  For a smooth function $u:M\longrightarrow R$, the \emph{gradient vector} of $u$ at $x$ is defined as the Legendre transform of the derivative of $u$, $\nabla u(x):=L^{\ast}(du(x))\in T_{x}M$.  Explicitly, we can write in coordinates
$$\nabla u(x):=\left\{\begin{array}{l}
                 g^{ij}(x,\nabla u)\frac{\partial u}{\partial x^{j}}\frac{\partial }{\partial x^{i}}~~~~~~ du(x) \neq 0,\\
                 0~~~~~~~~~~~~~~~~~~~~~~~~~~du(x)=0.
               \end{array}\right.\eqno(2.6)$$
It is $C^{\infty}$ on the open set $\{du\neq 0\}$ and $C^{0}$ at $\{du=0\}$.
Set $M_{V}:=\{x\in M|V(x)\neq 0\}$
for a vector field $V$ on $M$, and $M_{u}:=M_{\nabla u}$. For a
$C^{\infty}$ vector field $V$ on $M$ and $x\in M_{V}$, we define $\nabla V(x)\in T^{\ast}_{x}M\otimes T_{x}M$ by using the covariant derivative as
$$\nabla V(v):=D^{V}_{v}V(x)\in T_{x}M,~~~~~~~~v\in T_{x}M.\eqno(2.7)$$
We also set $\nabla^{2}u(x):=\nabla(\nabla u)(x)$ for the smooth function $u:M\longrightarrow R$ and $x\in M_{u}$. Let $\{e_{a}\}_{a=1}^{n}$ be a local orthonormal basis with respect to $g_{\nabla u}$ on $ M_{u}$. (In order to distinguish local orthonormal basis with respect to $g_{\nabla u}$ from that with respect to $g_{y}$, we use convention of the index range $1\leq a,~b,\cdots\leq n$.)  Using (2.5)-(2.7) and noting that $C_{\nabla u}(\nabla u,e_{a}, e_{b})=0$, we then have
$$\nabla^{2}u=\sum\left(\nabla^{2}u(e_{b})\right)\omega^{b}=\sum\left(D^{\nabla u}_{e_{b}}(\nabla u)\right)\omega^{b}$$
$$=\sum g_{\nabla u}\left(D^{\nabla u}_{e_{b}}(\nabla u),e_{a}\right)e_{a}\omega^{b}~~~~~~~$$
$$~~~~~~~~~~~~~~~~~~~=\sum\{e_{b}\left(g_{\nabla u}(\nabla u,e_{a})\right)-g_{\nabla u}(\nabla u,D^{\nabla u}_{e_{b}}e_{a})\}e_{a}\omega^{b}~~~$$
$$~~~~~~~~=\sum\{e_{b}(e_{a}(u))-\left(D^{\nabla u}_{e_{b}}e_{a}\right)(u)\}e_{a}\omega^{b}~~~~~~$$
$$=\sum u_{a|b}e_{a}\omega^{b},~~~~~~~~~~~~~~~~~~~~~~~~~$$
and
$$g_{\nabla u}(\nabla^{2}u(e_{a}),e_{b})=g_{\nabla u}(D^{\nabla u}_{e_{a}}(\nabla u),e_{b})~~~~~~~~~~~~~~~~~~~~~~~~~~~~~~~~~~~~~~~~~~~~~~~~$$
$$=e_{a}\left(g_{\nabla u}(\nabla u,e_{b})\right)-g_{\nabla u}(\nabla u,D^{\nabla u}_{e_{a}}e_{b})~~~$$
$$=e_{a}(e_{b}(u))-g_{\nabla u}(\nabla u,D^{\nabla u}_{e_{b}}e_{a}+[e_{a},e_{b}])$$
$$~~~~~~~~~~~~~~~~~~=e_{b}(e_{a}(u))+[e_{a},e_{b}](u)-g_{\nabla u}(\nabla u,D^{\nabla u}_{e_{b}}e_{a})-[e_{a},e_{b}](u)$$
$$=e_{b}\left(g_{\nabla u}(\nabla u,e_{a})\right)-g_{\nabla u}(\nabla u,D^{\nabla u}_{e_{b}}e_{a})~~~~$$
$$~~=g_{\nabla u}(D^{\nabla u}_{e_{b}}(\nabla u),e_{a})=g_{\nabla u}(\nabla^{2}u(e_{b}),e_{a}).~~$$
Namely,
$$u_{a|b}=u_{b|a},~~\forall a,b.$$

Next we define the \emph{divergence} of a $C^{\infty}$ vector field $V$ on $M$ with respect to an arbitrary volume form $d\mu$ by
$$ \textmd{div}V :=\sum_{i=1}^{n}(\frac{\partial V_{i}}{\partial x^{i}}+V_{i}\frac{\partial \Phi}{\partial x^{i}}),\eqno(2.8)$$
where $d\mu =e^{\Phi}dx^{1}dx^{2}\cdots dx^{n}$. Then the \emph{Finsler-Laplacian} of $u$ can be defined by
$$\Delta u:=\textmd{div}(\nabla u).\eqno(2.9)$$
Given a vector field $V$ such that $V\neq 0$ on $M_{u}$, we define the \emph{weighted gradient vector}$^{[13, 21]}$ and the  \emph{weighted Laplacian}$^{[13,21]}$ on the weighted Riemannian manifold $(M,g_{V})$ by
$$\nabla^{V}u:=\left\{\begin{array}{l}
                 g^{ij}(V)\frac{\partial u}{\partial x^{j}}\frac{\partial }{\partial x^{i}}~~~~~~on ~~~M_{u}, \\
                 0~~~~~~~~~~~~~~~~~~~~~on~~~M\backslash M_{u},
               \end{array}\right.
~~~~~~~~~~~~~~~\Delta^{V}u:=\textmd{div}(\nabla^{V}u).\eqno(2.10)$$
Clearly, the relation between the two gradients and that between the two Laplacians are
$$\nabla^{\nabla u}u=\nabla u,~~~~~~~~~~~~~~~~~\Delta^{\nabla u}u=\Delta u.$$

Let $(M,F,d\mu)$ be an $n-$dimensional Finsler manifold. If there is a constant $\lambda$ such that
$$\Delta f=-\lambda f$$
 for some function $f\in C^{2}(M)$, then the constant $\lambda$ is called the \emph{eigenvalue} of $\Delta$ and the function $f$ is called the \emph{eigenfunction} corresponding to $\lambda$. The least nonzero eigenvalue $\lambda_{1}$ of $\Delta$ is called the \emph{first eigenvalue} on $(M,F,d\mu)$. Let $\Omega\subset M$ be a domain with compact closure and nonempty boundary $\partial\Omega$. The first eigenvalue $\lambda_{1}(\Omega)$ of $\Omega$ is defined by$^{[15]}$
 $$\lambda_{1}(\Omega)=\inf_{u\in L^{1,2}_{0}(\Omega)}\frac{\int_{\Omega}(F^{\ast}(du))^{2}d\mu}{\int_{\Omega}u^{2}d\mu},$$
where $L^{1,2}_{0}(\Omega)$ is the completion of $C_{0}^{\infty}$ with respect to the norm
$$\|\varphi\|^{2}_{\Omega}=\int_{\Omega}\varphi^{2}d\mu+\int_{\Omega}(F^{\ast}(d\varphi))^{2}d\mu.$$
If $\Omega_{1}\subset\Omega_{2}$ are bounded domains, then $\lambda_{1}(\Omega_{1})\geq\lambda_{2}(\Omega_{2})\geq0$. Thus, if $\Omega_{1}\subset\Omega_{2}\subset\cdots\subset M$ are bounded domains so that $\bigcup\Omega_{i}=M$,  then the following limit
$$\lambda_{1}(M)=\lim_{i\rightarrow\infty}\lambda_{1}(\Omega_{i})\geq0$$
exists, and it is independent of the choice of $\{\Omega_{i}\}$.

In the end of this section, some lemmas are given below. \\

\textbf{\wuhao\heiti Lemma 2.1.(Bonnet-Myers)} \ \ \emph{Let $(M,F)$ be an $n-$dimensional forward geodesically complete connected Finsler manifold. If its Ricci curvature satisfies}
$$\textmd{Ric}\geq (n-1)k$$
\emph{for some positive constant $k$, then $M$ is compact and the diameter of $(M,F)$ is at most $\frac{\pi}{\sqrt{k}}$.}
\\

\textbf{\wuhao\heiti Lemma 2.2.([13])} \ \ \emph{Let $(M,F)$ be an $n-$dimensional Finsler manifold. Given $u\in C^{\infty}(M)$, we have}
$$\Delta^{\nabla u}(\frac{F(\nabla u)^{2}}{2})-D(\Delta u)(\nabla u)=\|\nabla u\|^{2}\textmd{Ric}_{\infty}(\nabla u)+\|\nabla^{2} u\|^{2}_{HS(\nabla u)}\eqno(2.11)$$
\emph{as well as}
$$\Delta^{\nabla u}(\frac{F(\nabla u)^{2}}{2})-D(\Delta u)(\nabla u)\geq \|\nabla u\|^{2}\textmd{Ric}_{N}(\nabla u)+\frac{(\Delta u)^{2}}{N}\eqno(2.12)$$
\emph{for $N\in[n,\infty]$, point-wise on $ M_{u}$. Here $\|\nabla^{2} u\|^{2}_{HS(\nabla u)}$ stands for the Hilbert-Schmidt norm with respect to $g_{\nabla u}$.}\\

According to Lemma 3.3 in [20], Lemma 3.2 in [13] and our discussion  on $\nabla^{2}u$ above, we can rewrite the result as\\

\textbf{\wuhao\heiti Lemma 2.3.} \ \ \emph{ Let $(M,F)$ be an $n-$dimensional Finsler manifold and $u:M\longrightarrow R$  a smooth function. Then on $M_{u}$ we have
$$\Delta u=\textmd{tr}_{g_{\nabla u}}(\nabla^{2}u)-S(\nabla u)=\sum_{a}u_{a|a}-S(\nabla u),\eqno(2.13)$$
where $u_{a|a}=g_{\nabla u}\left(\nabla^{2}u(e_{a}),e_{a}\right)$ and $\{e_{a}\}_{a=1}^{n}$ is a local $g_{\nabla u}$-orthonormal basis  on $ M_{u}$.}\\

\textbf{\wuhao\heiti Lemma 2.4.([17])} \ \ \ \emph{Let $(M,F,d\mu)$ be an $n-$dimensional complete connected Finsler manifold. Suppose that}
$$\textmd{Ric}\geq (n-1)k,~~~~~~~~~\|S\|\leq \Lambda.$$
\emph{Then for any $0<r<R$,}
$$\frac{\textmd{vol}^{d\mu}_{F}(B(x,R))}{V_{k,\Lambda,n}(R)}\leq\frac{\textmd{vol}^{d\mu}_{F}(B(x,r))}{V_{k,\Lambda,n}(r)},$$
\emph{where} $$ \|S\|_{x}:=\sup_{X\in T_{x}M\backslash0}\frac{S(X)}{F(X)};~~~~~~~~~~V_{k,\Lambda,n}(r):=\textmd{vol}(S^{n-1}(1))\int_{0}^{r}e^{\Lambda t}s_{k}(t)^{n-1}dt$$
 \emph{and $s_{k}$ denotes the unique solution to $y^{''}+ky=0$ with $y(0)=0,~y^{'}(0)=1.$}\\

\vspace*{-0.4cm}
\section*{\sihao\heiti 3.  \   Proofs of the main results}
 \textbf{\wuhao\heiti Theorem 3.1.} \ \ \ \emph{Let $(M,F)$ be an $n-$dimensional forward geodesically complete connected Finsler manifold. If the weighted Ricci curvature and $S$-curvature  satisfy}
$$\textmd{Ric}_{N} \geq (n-1)k,~~~~~~~~~~~~\dot{S} \leq \frac{(N-n)(n-1)}{N-1}k$$
\emph{for some uniform positive constant $k$ and $N\in (n,\infty)$, where $\dot{S}$ denotes the change rate of the $S$-curvature along
geodesics, then
$$\lambda_{1}\geq \frac{n-1}{N-1}Nk.$$
Moreover, the diameter of $M$ is $\sqrt{\frac{N-1}{n-1}}\frac{\pi}{\sqrt{k}}$ if the equality holds.}\\

\textbf{Proof}\ \ First of all, from (2.4) we see that, under the hypothesis in Theorem 3.1,
$$\textmd{Ric}=\textmd{Ric}_{N}-\dot{S}+\frac{S^{2}}{(N-n)F^{2}}\geq \frac{(n-1)^{2}}{N-1}k.\eqno(3.1)$$
So, $M$ is compact according to Lemma 2.1.

Let $u$ be the first eigenfunction on $(M,F)$ corresponding to the eigenvalue $\lambda_{1}$. This implies that
$$\Delta u=-\lambda_{1}u.$$
Furthermore, from the fact
$$\Delta^{\nabla u} u^{2}=\textmd{div} (\nabla^{\nabla u} u^{2})=\textmd{div}(2u\nabla u)=2u\Delta u+2\|\nabla u\|^{2}$$
we get
$$(\Delta u)^{2}=-\lambda_{1}u\Delta u=\lambda_{1}(\|\nabla u\|^{2}-\frac{1}{2}\Delta^{\nabla u} u^{2}).\eqno(3.2)$$
Integrating (2.12) and using divergence lemma on $M$, we obtain
$$\int_{M}\lambda_{1}\|\nabla u\|^{2}dV_{M}\geq\int_{M}\left(\|\nabla u\|^{2}\textmd{Ric}_{N}(\nabla u)+\frac{(\Delta u)^{2}}{N}\right)dV_{M}.$$
Thus, the assumption of the Theorem 3.1 and (3.2) yield
$$\int_{M}\left(\frac{N-1}{N}\lambda_{1}-(n-1)k\right)\|\nabla u\|^{2}dV_{M}\geq 0,$$
which means that
$$\lambda_{1}\geq\frac{n-1}{N-1}Nk.$$
If $\lambda_{1}=\frac{n-1}{N-1}Nk $, then all of the relevant inequalities become the equalities. We recall the formula (2.12), which was derived from (2.11) and the following inequalities.
$$\|\nabla^{2} u\|^{2}_{HS(\nabla u)}=\textmd{tr}(B(0)^{2})=\frac{(\textmd{tr}B(0))^{2}}{n}+\|B(0)-\frac{\textmd{tr}(B(0)}{n}I_{n}\|^{2}_{HS}~~~~$$
$$\geq\frac{(\textmd{tr}B(0))^{2}}{n}=\frac{\left(\Delta u+D\Psi(\nabla u)\right)^{2}}{n}~~~~~~~$$
$$~~~~~~~~~~~~~~~~~~~~~=\frac{(\Delta u)^{2}}{N}-\frac{\left(D\Psi(\nabla u)\right)^{2}}{N-n}+\frac{N(N-n)}{n}\left(\frac{\Delta u}{N}+\frac{D\Psi(\nabla u)}{N-n}\right)^{2}$$
$$\geq\frac{(\Delta u)^{2}}{N}-\frac{\left(D\Psi(\nabla u)\right)^{2}}{N-n},~~~~~~~~~~~~~~~~~\eqno(3.3)$$
where $B(0)=(\nabla^{2} u)\triangleq(u_{a|b})$ in the sense that $\nabla^{2} u (e_{a})=\sum _{b=1}^{n}u_{a|b}e_{b}$ (cf. [13], p.9,11-12),  $D\Psi(\nabla u)=S(\nabla u)$ by (2.4). So, under the condition $\lambda_{1}=\frac{n-1}{N-1}Nk $ we have
$$B(0)=\frac{tr(B(0)}{n}I_{n},\eqno(3.4)$$
$$\frac{\Delta u}{N}=-\frac{S(\nabla u)}{N-n}.\eqno(3.5)$$
Obviously  from (3.4) we can get
$$ u_{a|a}=u_{b|b},~~~~\forall a,b;~~~~~~~~~~~~         u_{a|b}=0,~~~~~~~\texttt{for}~~~a\neq b.\eqno(3.6)$$
Substituting (3.5) into (3.3), one has
$$\|\nabla^{2} u\|^{2}_{HS(\nabla u)}=\frac{\left(\Delta u+S(\nabla u)\right)^{2}}{n}=\frac{n}{N^{2}}(\Delta u)^{2}=\frac{n\lambda_{1}^{2}}{N^{2}}u^{2}.$$
Therefore combining (3.6) with the formula above, it holds that
$$ u_{a|a}^{2}=\frac{\lambda_{1}^{2}u^{2}}{N^{2}},~~\forall a.\eqno(3.7)$$
However, from Lemma 2.3  and (3.6) we also have
$$-\lambda_{1}u=\Delta u=nu_{a|a}-S(\nabla u),~~\forall a,$$
which together with (3.5) and (3.7) yields
$$u_{a|a}=-\frac{\lambda_{1} u}{N},~~\forall a.$$
Let $f(x)=\|\nabla u\|^{2}+\frac{\lambda_{1}}{N}u^{2}$. Then $f$ is $C^{\infty}$ on the open set $M_{u}$ and $C^{0}$ on $M\backslash M_{u}$. Its derivative  in the direction $e_{c}, \forall c$ on $M_{u}$ is
$$df(e_{c})=dg_{\nabla u}(\nabla u,\nabla u)(e_{c})+\frac{2\lambda_{1}}{N}uu_{c}$$
$$~~~~~~~~=2g_{\nabla u}(\nabla^{2} u,\nabla u)(e_{c})+\frac{2\lambda_{1}}{N}uu_{c}$$
$$~~~~~~~~~~~~~~~~~~~~~~=2g_{\nabla u}(\sum u_{a|b}e_{a}\omega^{b},\sum u_{d}e_{d})(e_{c})+\frac{2\lambda_{1}}{N}uu_{c}$$
$$~~~~~~~~~~~~~~~~~~~~=2u_{c}u_{c|c}+\frac{2\lambda_{1}}{N}uu_{c}=0,~~~~~~~~~~~~~~~~~~~~~~$$
which means that $f$ is constant on $ M_{u}$. On the other hand, we claim that $f$ is also constant on $M\backslash M_{u}$. In fact, if $M\backslash M_{u}\ni x$ is an inner point, then $f=\frac{\lambda_{1}}{N}u^{2}$ holds on a neighborhood $U$ of $x$ so that $df=0$ or $f$ is constant on $U$. If $M\backslash M_{u}\ni x$ is a boundary point,  we choose a sequence $\{x_{k}\}\subset M_{u}$ such that $x_{k}\longrightarrow x, (k\longrightarrow\infty)$. Then $f(x)=f|_{M_{u}}$ according to the continuity of $f$. Finally, using the continuity of $f$ again and connectivity of $M$ we obtain that the  function $f(x)$ is constant on $M$.

Suppose that $u$ attains its maximum $u_{max}$ and minimum $u_{min}$ at $p\in M$ and $q\in M$ respectively. Since $\|\nabla u\|^{2}=0$ at both $p$ and $q$, we see that $f(p)=\frac{\lambda_{1}}{N}(u_{max})^{2}=f(q)=\frac{\lambda_{1}}{N}(u_{min})^{2}$, which implies that $|u_{max}|=|u_{min}|$. This also mean that all maximum (or minimum) of $u$ are equal. Without loss of generality, we can assume that $u_{max}=1$ and $u_{min}=-1$. Let $\gamma(s)$ be the minimal regular geodesic of $(M,F)$ from $p$ to $q$ with the tangent vector $\dot{\gamma}(s)$. We can suppose that along $\gamma(s)$ there is not any other extreme point. Otherwise, Since $u$ is continuous, $p$ must not be the cluster point of minimal extreme points of $u$. Hence we may assume $q^{'}\in\gamma(s)$ is the first minimal extreme point of $u$ from $p$. Next set off from $q^{'}$ to $p$ along $\overleftarrow{\gamma(s)}$, by the same way we get the maximum  extreme point $p^{'}\in\gamma(s)$. Then $\gamma(s)|_{\widehat{p^{'}q^{'}}}$ is the minimal regular geodesic without other extreme point of $u$. So we might as well assume that $\gamma(s)$ has this property which means $\|\nabla u\|(x)>0, \forall x\in\gamma(s)\backslash(p,q)$.  Consequently $\gamma(s)\backslash \{p,q\}\subset M_{u}$.
 Since  $\lambda_{1}=\frac{n-1}{N-1}Nk $, then we have $\frac{\|\nabla u\|}{\sqrt{1-u^{2}}}=\sqrt{\frac{n-1}{N-1}k}$ along $\gamma(s)$.

Let $d_{M}$ denote the diameter of $(M,F)$. We then have
$$\sqrt{\frac{n-1}{N-1}k}d_{M}\geq\sqrt{\frac{n-1}{N-1}k}\int_{\gamma}F(\dot{\gamma})ds=\int_{\gamma}F(\dot{\gamma})\frac{\|\nabla u\|}{\sqrt{1-u^{2}}}ds.\eqno(3.8)$$
From $|\frac{du}{ds}|=|g_{\nabla u}(\nabla u, \dot{\gamma})|\leq F(\dot{\gamma})\|\nabla u\|$ one gets
$$\int_{\gamma}F(\dot{\gamma})\frac{\|\nabla u\|}{\sqrt{1-u^{2}}}ds\geq\int_{-1}^{1}\frac{du}{\sqrt{1-u^{2}}}=\pi.\eqno(3.9)$$
It follows from (3.8) and (3.9) that $d_{M}\geq\sqrt{\frac{N-1}{n-1}}\frac{\pi}{\sqrt{k}}$.

On the other hand, from (3.1) and Lemma 2.1 we can obtain $d_{M}\leq\sqrt{\frac{N-1}{n-1}}\frac{\pi}{\sqrt{k}}$. So $d_{M}=\sqrt{\frac{N-1}{n-1}}\frac{\pi}{\sqrt{k}}$. This finish the proof.\\

From the proof of Theorem 3.1, it is not difficult to obtain\\

 \textbf{\wuhao\heiti Proposition 3.2.} \ \ \ \emph{Let $(M,F)$ be an $n-$dimensional compact connected Finsler manifold. If the weighted Ricci curvature satisfies} $\textmd{Ric}_{N} \geq (n-1)k$
\emph{for some uniform positive constant $k$ and $N\in (n,\infty)$, then
$$\lambda_{1}\geq \frac{n-1}{N-1}Nk.$$
Moreover, the diameter of $M$ is at least $\sqrt{\frac{N-1}{n-1}}\frac{\pi}{\sqrt{k}}$ if the equality holds.}\\

\textbf{\wuhao\heiti Theorem 3.3.} \ \ \ \emph{Let $(M,F)$ be an $n-$dimensional forward geodesically complete connected Finsler manifold. If $S=0$ and Ricci curvature }
$\textmd{Ric} \geq (n-1)k$
\emph{for some uniform positive constant $k$, then
$$\lambda_{1}\geq nk.$$
Moreover, if the equality holds, then the diameter of $M$ is $\frac{\pi}{\sqrt{k}}$,  and $M$ is homeomorphic to $S^{n}$. In particular, if $F$
is reversible and $M$ has Busemann-Hausdorff volume form, then $(M,F)$ is isometric to $S^{n}(\frac{1}{\sqrt{k}})$.}\\

\textbf{Proof}\ \ If $S=0$, then $\textmd{Ric}_{N}=\textmd{Ric}$ from the Definition 2.1. Therefore, by Theorem 3.1 we can easily get the first part of Theorem 3.3.  Next we only prove the last part when the equality holds. Under the condition of Theorem 3.2, $f(x)=\|\nabla u\|^{2}+\frac{\lambda_{1}}{n}u^{2}$. Here $f(x)$ is constant on $M$ by the proof of Theorem 3.1.  Put
$$M^{+}=\{x\in M|u(x)>0\},~~~M^{0}=\{x\in M|u(x)=0\},~~~M^{-}=\{x\in M|u(x)<0\}.$$
Then $M^{+},~M^{-}$ are open sets on $M$, and $M^{0}$ is a close set with zero measure. Let $p$ and $q$ are the maximal point and minimal point of $u$ respectively with $u(p)=1,~u(q)=-1$. So, if $\lambda_{1}= nk,$ then $\frac{\|\nabla u\|}{\sqrt{1-u^{2}}}=\sqrt{k}$. Suppose that $\gamma$ is the minimal geodesic  of $(M,F)$ from $p$ to $q$ with the tangent vector $\dot{\gamma}(s)$. Denote by $L(\gamma)$ the length of  $\gamma$. Then
$$\sqrt{k}L(\gamma)=\int_{\gamma}F(\dot{\gamma})\frac{\|\nabla u\|}{\sqrt{1-u^{2}}}ds\geq\int_{-1}^{1}\frac{du}{\sqrt{1-u^{2}}}=\pi\eqno(3.10)$$
which means that $L(\gamma)=d(p,q)=d$. Similarly, we also get $d(q,p)=d$. Furthermore, we claim $B(p,\frac{d}{2})\subset M^{+}$. In fact, if there exists a point $x_{0}\in M^{-}\cup M^{0}$ such that $x_{0}\in B(p,\frac{d}{2})$, then we  suppose that $\eta$ is the minimal geodesic  of $(M,F)$ from $p$ to $x_{0}$ with the tangent vector $\dot{\eta}(s)$. Thus
$$\sqrt{k}L(\eta)=\int_{\eta}F(\dot{\eta})\frac{\|\nabla u\|}{\sqrt{1-u^{2}}}ds\geq\int_{0}^{1}\frac{du}{\sqrt{1-u^{2}}}=\frac{\pi}{2}\eqno(3.11)$$
which shows that $L(\eta)=d(p,x_{0})\geq\frac{d}{2}$. This contradict the assumption. Similarly, $B(q,\frac{d}{2})\subset M^{-}$. So we get
 $$B(p,\frac{d}{2})\cap B(q,\frac{d}{2})=\emptyset.\eqno(3.12)$$
Note that if $S=0, k>0$, then $V_{k,\Lambda,n}(r)=\textmd{vol}(S^{n}(k;r))$. Hence from Lemma 2.4 we get
$$\frac{\textmd{vol}^{d\mu}_{F}(B(p,\frac{\pi}{2\sqrt{k}}))}{\textmd{vol}(S^{n}(k;\frac{\pi}{2\sqrt{k}}))}\geq \frac{\textmd{vol}^{d\mu}_{F}(B(p,\frac{\pi}{\sqrt{k}}))}{\textmd{vol}(S^{n}(k;\frac{\pi}{\sqrt{k}}))}
=\frac{\textmd{vol}^{d\mu}_{F}M}{\textmd{vol}S^{n}(\frac{1}{\sqrt{k}})},$$
which implies that
$$\textmd{vol}^{d\mu}_{F}(B(p,\frac{d}{2}))\geq\frac{1}{2}\textmd{vol}^{d\mu}_{F}M.\eqno(3.13)$$
A similar argument yields
$$\textmd{vol}^{d\mu}_{F}(B(q,\frac{d}{2}))\geq\frac{1}{2}\textmd{vol}^{d\mu}_{F}M.\eqno(3.14)$$
From (3.12), (3.13) and (3.14), we have
$$B(p,\frac{d}{2})= M^{+},~~~~~~~~~~~B(q,\frac{d}{2})= M^{-},\eqno(3.15)$$
 $ M^{0}$ is the boundary of both $B(p,\frac{d}{2})$ and $B(q,\frac{d}{2})$. In addition, we can prove that for any point $x\in M^{0}$, $d(p,x)=\frac{d}{2}$. On the one hand, from (3.11), $d(p,x)\geq\frac{d}{2}$. On the other hand, if $d(p,x)>\frac{d}{2}$, then there exists a neighborhood $U$ of $x$ such that $d(p,y)>\frac{d}{2}$ for any $y\in U$. This contradict (3.15). Similarly, for any point $x\in M^{0}$, $d(q,x)=\frac{d}{2}$.

In the following, we illustrate that $u$ has only one maximal point on $M$. If not, we assume $p_{1},p_{2}$ are the two maximal points of $u$.
Let $\sigma_{1}$ be the minimal regular geodesic from  $p_{1}$ to $q$. Set $x_{1}= \sigma_{1}\cap M^{0}$, then $L(\sigma_{1})=d(p_{1},q)=d$, $d(p_{1},x_{1})=L(\sigma_{1}|_{\widehat{p_{1}x_{1}}})=d(x_{1},q)=L(\sigma_{1}|_{\widehat{x_{1}q}})=\frac{d}{2}$. Draw a  minimal regular geodesic $\eta$ from  $p_{2}$ to $x_{1}$. Then $d(p_{2},x_{1})=L(\eta)=\frac{d}{2}$. From (3.10) we have
$$d(p_{2},x_{1})+d(x_{1},q)=d(p_{1},q).$$
Let $\sigma_{2}\triangleq \eta\cup\sigma_{1}|_{\widehat{x_{1}q}}$, then $\sigma_{2}$ is a minimal regular geodesic from $p_{2}$ to $q$ with $L(\sigma_{2})=d$. Note that the equality  in (3.10) holds if and only if $\dot{\gamma}$ is parallel to $\nabla u$ and $u$ is monotone decreasing along $\gamma$. Hence at $x_{1}$, we have $\dot{\sigma}_{1}(x_{1})=\dot{\sigma}_{2}(x_{1})=-\frac{\nabla u}{\|\nabla u\|}(x_{1})$. According to the uniqueness of geodesic we have $\sigma_{1}=\sigma_{2}$ so that $p_{1}=p_{2}$. Similarly, $u$ has only one minimal point $q$ on $M$.

Since $\|\nabla u\|^{2}+ku^{2}=k$, then we have
$$D^{\nabla u}_{\nabla u}\left(\frac{\nabla u}{\|\nabla u\|}\right)=D^{\nabla u}_{\nabla u}\left(\frac{\nabla u}{\sqrt{k}\sqrt{1-u^{2}}}\right)
~~~~~~~~~~~~~~~~~~~~~~~~~~~~~~~~~~~~~$$
$$~~~~~~~~~~~~=\frac{1}{\sqrt{k}\sqrt{1-u^{2}}}D _{\nabla u}^{\nabla u}\nabla u+D _{\nabla u}^{\nabla u}\left(\frac{1}{\sqrt{1-u^{2}}}\right)\frac{\nabla u}{\sqrt{k}}$$
$$~~~~~~~~~~~~~~~~~~~~~~~~~~~=\frac{1}{\sqrt{k}\sqrt{1-u^{2}}}\nabla^{2}u(\nabla u)+g_{\nabla u}\left(\nabla u,\nabla^{\nabla u}\left(\frac{1}{\sqrt{1-u^{2}}}\right)\right)\frac{\nabla u}{\sqrt{k}}$$
$$=\frac{1}{\sqrt{k}\sqrt{1-u^{2}}}(\nabla^{2}u(\nabla u)+uk\nabla u)=0,$$
which means that $\frac{\nabla u}{\|\nabla u\|}$ is geodesic field. For any $x_{0}\in M$, Draw a minimal geodesic $\gamma$ from $q$ to $x_{0}$, then
 $$\sqrt{k}L(\gamma)=\int_{\gamma}F(\dot{\gamma})\frac{\|\nabla u\|}{\sqrt{1-u^{2}}}ds\geq\int_{-1}^{u(x_{0})}\frac{du}{\sqrt{1-u^{2}}}.$$
Since  $\gamma$ is minimal geodesic, then $\dot{\gamma}=\frac{\nabla u}{\|\nabla u\|}$. Further, we have on $\gamma$
$$|u^{'}|^{2}+ku^{2}=k,~~~~~  u(0)=-1,~~u^{'}(0)=0,$$
which shows that $u=-cos\sqrt{k}t,t\in[0,\frac{\pi-\arccos u(x_{0})}{\sqrt{k}}]$. As a geodesic on $M$, $\gamma$ is defined in $[0,\infty]$, so we have  $u=-cos\sqrt{k}t,t\in[0,\frac{\pi}{\sqrt{k}}]$. Particularly, $u(\gamma(\frac{\pi}{\sqrt{k}}))=1$ which means $p\in\gamma$. Clearly, the point $p$ is the cut
locus of $q$. Thus we conclude that $\exp_{q}:T_{q}M\supset B_{q}(\frac{\pi}{\sqrt{k}})\longrightarrow M^{n}\backslash\{p\}$ is diffeomorphism. On the other hand, $\exp_{\tilde{q}}:T_{\tilde{q}}S^{n}\supset B_{\tilde{q}}(\pi)\longrightarrow S^{n}\backslash\{\tilde{p}\}$ is also diffeomorphism where $S^{n}$ is $n$-sphere, $\tilde{q},\tilde{p}$ are the south pole and north pole respectively. Let $(\tilde{r},\tilde{\theta}^{\alpha})$ be the polar coordinate system of $T_{\tilde{q}}S^{n}$ and $(r,\theta^{\alpha})$ be  the polar coordinate system of $T_{q}M^{n}$. Define $h:T_{\tilde{q}}S^{n}\longrightarrow T_{q}M$ by $r=\frac{\tilde{r}}{\sqrt{k}},\theta^{\alpha}=\tilde{\theta}^{\alpha}$, then $h$ is diffeomorphism. Now we define $\psi:M^{n}\longrightarrow S^{n}$ by
$$\psi(x)=\left\{\begin{array}{cc}
            \exp_{\tilde{q}}\circ h^{-1}\circ \exp_{q}^{-1}(x) & x\neq p \\
            \tilde{p} & x=p
          \end{array}\right.$$

It is not hard to see $\psi$ is  homeomorphic. i.e. $M$ is homeomorphic to $S^{n}$. At last,  if  $F$ is reversible, $S_{BH}=0$ and
   the diameter of $M$ is $\frac{\pi}{\sqrt{k}}$, then according to the Corollary 1 in [8], $(M,F)$ is isometric to $S^{n}(\frac{1}{\sqrt{k}})$.
    The theorem has been proved.\\

\textbf{\wuhao\heiti Theorem 3.4.} \ \ \  \emph{Let $(M,F)$ be an $n-$dimensional compact Finsler manifold. If the weighted Ricci curvature } $\textmd{Ric}_{\infty} \geq 0$,
\emph{then
$$\lambda_{1}\geq\frac{\pi^{2}}{d^{2}},$$
where $d$ denotes the diameter of $(M,F)$.}\\

\textbf{Proof}\ \ Let $u$ be the first eigenfunction on $(M,F)$ corresponding to the first eigenvalue $\lambda_{1}$. Since $\int_{M}ud\mu=-\frac{1}{\lambda_{1}}\int_{M}\Delta ud\mu=0$ and noting that $-u$ is not necessarily  the first eigenfunction on $(M,F)$, we have to
 assume that there are two case: $\sup u=1$ or $\inf u=-1$.

\textbf{Case (I):}
$$1=\sup u>\inf u=-k\geq-1,~~~~~~~~~~~~0<k\leq1.$$
For small $\varepsilon>0$, let
$$v=\frac{u-\frac{1}{2}(1-k)}{\frac{1}{2}(1+k)(1+\varepsilon)}.$$
Clearly, $dv=\frac{2}{(1+k)(1+\varepsilon)}du$. Since  Legendre transform $L^{\ast}:T^{\ast}M\longrightarrow TM$ is dimorphism and satisfies $L^{\ast}(a\zeta)=aL^{\ast}(\zeta), a\in R^{+}, \zeta\in T^{\ast}M$, we have $$\nabla v=\nabla^{\nabla u} v=\frac{2}{(1+k)(1+\varepsilon)}\nabla u$$
under which
$$\left\{\begin{array}{c}
    \Delta v=-\lambda_{1}(v+a_{\varepsilon}),~~~~~a_{\varepsilon}=\frac{1-k}{(1+k)(1+\varepsilon)},\\
    \sup v=\frac{1}{1+\varepsilon} ,~~~~~~~~\inf v=-\frac{1}{1+\varepsilon}.
  \end{array}\right.$$
  
\textbf{Case (II):}
$$1\geq k=\sup u>\inf u=-1,~~~~~~~~~~~~0<k\leq1.$$
For small $\varepsilon>0$, let
$$v=\frac{u+\frac{1}{2}(1-k)}{\frac{1}{2}(1+k)(1+\varepsilon)}.$$
Then we also have $$\nabla v=\nabla^{\nabla u} v=\frac{2}{(1+k)(1+\varepsilon)}\nabla u,$$
under which
$$\left\{\begin{array}{c}
    \Delta v=-\lambda_{1}(v-a_{\varepsilon}),~~~~~a_{\varepsilon}=\frac{1-k}{(1+k)(1+\varepsilon)},\\
    \sup v=\frac{1}{1+\varepsilon} ,~~~~~~~~\inf v=-\frac{1}{1+\varepsilon}.
  \end{array}\right.$$
Let $v=\sin \theta$, then
$$-\frac{1}{1+\varepsilon}\leq\sin \theta\leq\frac{1}{1+\varepsilon},~~~~~~\frac{\|\nabla v\|^{2}}{1-v^{2}}=\|\nabla^{\nabla u} \theta\|^{2}.$$
Consider the function
$$f(x)=\frac{\|\nabla v\|^{2}}{1-v^{2}}.$$

Since $M$ is compact, then we can apply the maximal principle to $f(x)$ on the weighted Riemannian manifold $(M,g_{\nabla u})$. Suppose that $f(x)$ attains its maximum at $x_{0}\in M$, then $\nabla^{\nabla u}f(x_{0})=0$, $\Delta^{\nabla u}f(x_{0})\leq0$ and $x_{0}\in M_{u}$.

Let $\{e_{a}\}_{a=1}^{n}$ be a local orthonormal basis with respect to $g_{\nabla u}$ on $ M_{u}$.   Write
 $\nabla v=\displaystyle{\sum_{a}} v_{a}e_{a}$. Then by simple computations on $\nabla^{\nabla u}f(x_{0})=0$ we have
$$\sum_{b}v_{b}v_{b|a}=\frac{\|\nabla v\|^{2}(-v)v_{a}}{1-v^{2}},~~~~\forall a.\eqno(3.16)$$
Furthermore, a straightforward calculation yields
$$\Delta^{\nabla u}f(x_{0})=\Delta^{\nabla u}\left(\frac{\|\nabla v\|^{2}}{1-v^{2}}\right)=\frac{\Delta^{\nabla u}\left(\|\nabla v\|^{2}\right)}{1-v^{2}}+\|\nabla v\|^{2}\Delta^{\nabla u}\left(\frac{1}{1-v^{2}}\right)$$
$$~~~~~~~~~~~~~~~~~+2g_{\nabla u}\left(\nabla^{\nabla u}\left(\|\nabla v\|^{2}\right),\nabla^{\nabla u}\left(\frac{1}{1-v^{2}}\right)\right)\triangleq A+B+C,\eqno(3.17)$$
where
$$A=\frac{\Delta^{\nabla u}\left(\|\nabla v\|^{2}\right)}{1-v^{2}},~~~~~~~~~~~~~~~~~~~~~~~~~~~~~~~~~~~~~~~~~~~~~~~$$
$$B=\|\nabla v\|^{2}\Delta^{\nabla u}\left(\frac{1}{1-v^{2}}\right)=\|\nabla v\|^{2}\textmd{div}\left( \nabla^{\nabla u}\left(\frac{1}{1-v^{2}}\right)\right)$$
$$~~~~~~~~~~~~~~~~~~~~=\|\nabla v\|^{2}\left\{\frac{2v}{(1-v^{2})^{2}}\textmd{div}\left(\nabla^{\nabla u}v\right)+2g_{\nabla u}\left(\nabla^{\nabla u}v,\nabla^{\nabla u}\left(\frac{v}{(1-v^{2})^{2}}\right)\right)\right\}$$
$$=\|\nabla v\|^{2}\left\{\frac{2v}{(1-v^{2})^{2}}\Delta v+\frac{2\|\nabla v\|^{2}}{(1-v^{2})^{2}}+\frac{8v^{2}\|\nabla v\|^{2}}{(1-v^{2})^{3}}\right\},~~$$
$$C=2g_{\nabla u}\left(\nabla^{\nabla u}\left(\|\nabla v\|^{2}\right),\nabla^{\nabla u}\left(\frac{1}{1-v^{2}}\right)\right)=\frac{8vv_{a}v_{b}v_{a|b}}{(1-v^{2})^{2}}.$$
So, from the formulas above, we can rewrite (3.17) as follows
$$0\geq\Delta^{\nabla u}f(x_{0})=\frac{\Delta^{\nabla u}(\|\nabla v\|^{2})}{1-v^{2}}+\frac{8v\sum v_{a}v_{b}v_{a|b}}{(1-v^{2})^{2}}~~~~~~~~~~$$
$$~~~~~~~~~~~~~~~~~~~~~~~~~~~~-\frac{2\|\nabla v\|^{4}+2v\|\nabla v\|^{2}\Delta v}{(1-v^{2})^{2}}+\frac{8v^{2}\|\nabla v\|^{4}}{(1-v^{2})^{3}}.$$
Substituting (3.16) into it, one has
$$0\geq \Delta^{\nabla u}(\|\nabla v\|^{2})+\frac{2\|\nabla v\|^{4}+2v\|\nabla v\|^{2}\Delta v}{1-v^{2}}.\eqno(3.18)$$
From (2.11) and the conditions of Theorem 3.4, we get
$$~~~~~~\Delta^{\nabla u}(\|\nabla v\|^{2})=2\|\nabla v\|^{2}\textmd{Ric}_{\infty}(\nabla v)+2D(\Delta v)(\nabla v)+2\|\nabla^{2} v\|^{2}_{HS(\nabla v)}$$
$$~~\geq 2D(-\lambda_{1}(v\pm a_{\varepsilon}))(\nabla v)+2\sum_{ab} v_{a|b}^{2}$$
$$=-2\lambda_{1}\|\nabla v\|^{2}++2\sum_{ab} v_{a|b}^{2}.~~~~~~~~~\eqno(3.19)$$
By the Schwartz inequality and (3.16), we have
$$\sum_{ab} v_{a|b}^{2}\sum_{b} v_{b}^{2}\geq\sum_{a}(\sum_{b}v_{b}v_{b|a})^{2}=\sum_{a}\frac{\|\nabla v\|^{4}v^{2}v_{a}^{2}}{(1-v^{2})^{2}},$$
which means that
$$\sum_{ab} v_{a|b}^{2}\geq\frac{\|\nabla v\|^{4}v^{2}}{(1-v^{2})^{2}}.\eqno(3.20)$$
Utilizing (3.18)-(3.20) above, we obtain  at the point $x_{0}$ in both cases that
$$f(x_{0})=\frac{\|\nabla v\|^{2}}{1-v^{2}}(x_{0})\leq\lambda_{1}(1+a_{\varepsilon}).$$
So for any $x\in M$, we have
$$\sqrt{f(x)}=\|\nabla^{\nabla u} \theta\|\leq\sqrt{\lambda_{1}(1+a_{\varepsilon})}.\eqno(3.21)$$
Set
$$G(\theta)=\max _{\begin{array}{l}
                    x\in M \vspace*{-0.2cm}\\
                     \theta(x)=\theta
                   \end{array}}
\|\nabla^{\nabla u} \theta\|^{2}=\max _{\begin{array}{l}
                    x\in M \vspace*{-0.2cm}\\
                     \theta(x)=\theta
                   \end{array}}\frac{\|\nabla v\|^{2}}{1-v^{2}}.$$
Clearly, $G(\theta)\in C^{0}\left([-\frac{\pi}{2}+\delta,\frac{\pi}{2}-\delta]\right)$, where $\delta$ is specified by
$$\sin \left(\frac{\pi}{2}-\delta\right)=\frac{1}{1+\varepsilon},~~~~~G\left(-\frac{\pi}{2}+\delta\right)=G\left(\frac{\pi}{2}-\delta\right)=0.$$
From (3.21) we can write
$$G(\theta)\leq\lambda_{1}(1+a_{\varepsilon}).$$
Under which we let
$$G(\theta)=\lambda_{1}(1+a_{\varepsilon}\varphi(\theta)),~~~~~~\varphi(\theta)\in C^{0}\left([-\frac{\pi}{2}+\delta,\frac{\pi}{2}-\delta]\right).$$
Since  $G(\theta)$ vanishes at the end points of the interval $[-\frac{\pi}{2}+\delta,\frac{\pi}{2}-\delta]$, then
$$\varphi(\frac{\pi}{2}-\delta)=\varphi(-\frac{\pi}{2}+\delta)<-1.$$
By (3.21) we see that $\varphi(\theta)\leq1.$

In the following,  by the same way in [22], we can  get
$$\varphi(\theta)\leq\psi(\theta),\eqno(3.22)$$
where $\psi(\theta)$ is defined by
$$\psi(\theta)=\left\{\begin{array}{c}
                 \frac{\frac{4}{\pi}(\theta+\cos\theta\sin\theta)-2\sin\theta}{\cos^{2}\theta},~~~~~~~~~~\theta\in\left(-\frac{\pi}{2},\frac{\pi}{2}\right) \\
                 \psi(\frac{\pi}{2})=1,~~~~~~~~\psi(-\frac{\pi}{2})=-1.
               \end{array}\right.\eqno(3.23)$$
Now we continue to prove Theorem 3.2. From (3.22), we have
$$\|\nabla^{\nabla u} \theta\|\leq \sqrt{\lambda_{1}}\sqrt{1+a_{\varepsilon}\psi(\theta)}.\eqno(3.24)$$
Let $p, q \in M$ be such points that $\theta(p)=-\frac{\pi}{2}+\delta,~~\theta(q)=\frac{\pi}{2}-\delta$. Let $\gamma$ be a shortest geodesic joining $p$ and $q$.  Denote by $T$ the tangent vector of $\gamma$. Then
$$\|\nabla^{\nabla u} \theta\|=\frac{\|\nabla v\|}{\cos\theta}=\frac{F(\nabla v)}{\cos\theta}\geq\frac{|g_{\nabla u}\left(\nabla v,\frac{T}{F(T)}\right)|}{\cos\theta}$$
$$~~~~~~~~=\frac{|Tv|}{F(T)\cos\theta}=\frac{|\frac{dv}{ds}|}{F(T)\cos\theta}=\frac{\frac{d\theta}{ds}}{F(T)}.\eqno(3.25)$$
 Therefore from (3.24) and (3.25) one gets
$$\sqrt{\lambda_{1}}d\geq\int_{\gamma}\sqrt{\lambda_{1}}F(T)ds\geq\int_{-\frac{\pi}{2}+\delta}^{\frac{\pi}{2}
-\delta}\frac{d\theta}{\sqrt{1+a_{\varepsilon}\psi(\theta)}}.\eqno(3.26)$$
 It is easy to see from (3.23) that $\psi(0)=0, \psi(-\theta)=-\psi(\theta),| a_{\varepsilon}\psi(\theta)|<1 $. Hence, we have
 $$\int_{-\frac{\pi}{2}+\delta}^{\frac{\pi}{2}-\delta}\frac{d\theta}{\sqrt{1+a_{\varepsilon}\psi(\theta)}}
 =\int_{0}^{\frac{\pi}{2}-\delta}\left(\frac{1}{\sqrt{1+a_{\varepsilon}\psi(\theta)}}+\frac{1}{\sqrt{1-a_{\varepsilon}\psi(\theta)}}\right)d\theta~~~$$
$$~~~~~~~~~~~~~~~~~~~~~~~~~=2\int_{0}^{\frac{\pi}{2}-\delta}\left(1+\sum_{i=1}^{\infty}\frac{1\cdot3\cdots(4i-1)}
{2\cdot4\cdots4i}a_{\varepsilon}^{2i}\psi^{2i}\right)d\theta$$
$$\geq2(\frac{\pi}{2}-\delta)=\pi-2\delta.~~~~~~~~~\eqno(3.27)$$
Thus
$$\sqrt{\lambda_{1}}d\geq\pi-2\delta.$$
Letting $\varepsilon\rightarrow0$, so that $\delta\rightarrow0$ too, we then obtain
$$\lambda_{1}\geq\frac{\pi^{2}}{d^{2}}.$$

\textbf{Remark}\ \ The estimate in Theorem 3.4 has been pointed out in [21], where a sharp lower
bound for first Neumann eignenvalue of Finsler-Laplacian was given. The conclusion of Theorem
3.4 is not sharp for $n\geq2$.\\

 If  $S=(n+1)cF$ for some constant $c$, then $\dot{S}=0$ so that $\textmd{Ric}_{\infty}=\textmd{Ric}$. So we can easily get the following\\

\textbf{\wuhao\heiti Corollary 3.5.} \ \ \  \emph{Let $(M,F)$ be an $n-$dimensional compact Finsler manifold. If $M$ has constant $S$-curvature and } $\textmd{Ric}\geq 0$,
\emph{then
$$\lambda_{1}\geq\frac{\pi^{2}}{d^{2}},$$
where $d$ denotes the diameter of $(M,F)$.}\\

\textbf{\wuhao\heiti Acknowledgments}\\
 This project is supported by the National Natural Science Foundation of China under grant numbers 10971239, 11171253 and the Natural Science Foundation of High Education in Anhui Province under grant number KJ2012B197.

\vspace*{0.5cm}
\setlength{\parindent}{1em}
\def\hang{\hangindent\parindent}
\def\textindent#1{\noindent\llap{#1\enspace}\ignorespaces}
\def\re{\hang\textindent}
\vspace*{-0.3cm}

\begin{flushleft}
 \textbf{{\wuhao\heiti References\\ }}
\end{flushleft}
\bigskip
\setlength{\baselineskip}{1.2em} \vspace*{-0.6cm}
\xiaowuhao\songti

\re{}[1]\ P. Antonelli and B. Lacky, \emph{The theory of Finslerian Laplacian and application}. Math. and its appl. 459, Kluwer academic publishers, 1998.

\re{}[2]\ D. Bao, S. S. Chern and Z. Shen,  \emph{An introduction to Riemann-Finsler geometry}. GTM, 200. Springer, Berlin Heidelberg New York, 2000.

\re{}[3]\  D. Bao and B. Lacky, \emph{A Hodge decomposition theorem for Finsler spaces}. C. R. Acad. Sc. Paris 223(1996), 51-56.

\re{}[4]\ P. Centore, \emph{Finsler Laplacians and minimal-energy map}. Inter. J. Math. 11(2000), 1-13.

\re{}[5]\ I. Chavel, \emph{Eigenvalues in Riemannian geometry}.  Acad. Press, Inc., London, 1984.

\re{}[6]\  Y. Ge and Z. Shen, \emph{Eigenvalues and eigenfunctions of metric measure manifolds}. Proc. London Math. Soc.  82(2001), 725-746.

\re{}[7] F. Hang and X.Wang, A remark on Zhong-Yang's eigenvalue estimate. Int. Math. Res. Not. 18(2007), Art. ID rnm064, 9pp.

\re{}[8]\ C.-W. Kim and J.-W. Yim, \emph{Finsler manifolds with positive constant flag curvature}. Geom. Dedicata, 98(2003), 47-56.

\re{}[9]\  P. Li and S-Y. Yau, \emph{Estimates of eigenvalues of a compact Riemannian manifold, Geometry of the
Laplace operator} (Proc. Sympos. Pure Math., Univ. Hawaii, Honolulu, Hawaii, 1979), Proc. Sympos.
Pure Math., XXXVI, Amer. Math. Soc., Providence, R.I. 1980, 205-239.

\re{}[10]\   A. Lichnerowicz, \emph{Geometrie des groupes de transforamtions. Travaux et Recherches Mathemtiques}.
III. Dunod, Paris, 1958.

\re{}[11]\ J. Lott and C.Villani, \emph{Ricci curvature for metric-measure spaces via optimal transport}. Ann. of Math. 169(2009), 903-991.

\re{}[12]\ M. Obata, \emph{Certain conditions for a Riemannian manifold tobe isometric with a sphere}. J. Math. Soc. Japan. 14(1962), 333-340.

\re{}[13]\  S. Ohta and K-T, Sturm,  \emph{Bochner-Weitzenbock formula and Li-Yau estimates on Finsler manifolds}. arXiv: 1105.0983.

\re{}[14]\ S. Ohta, \emph{Finsler interpolation inequalities}. Calc. Var. Partial Differential Equations 36(2009), 211-249.

\re{}[15]\ Z. Shen, \emph{Lectures on Finsler geometry}. World Scientific Publishing Co., Singapore, 2001.

\re{}[16]\ ---------, \emph{The non-linear Laplacian for Finsler manifolds, "The theory of Finslerian Laplacians and applications}"  (edited by P.Antonelli), Proc. Conf. \emph{On Finsler Laplacians}, Kluwer Acad. Press, Netherlands, 1998.

\re{}[17]\ ---------, \emph{Volume comparison and its applications in Riemann-Finsler geometry}. Adv. Math. 128(1997), 306-328.

\re{}[18]\  K.-T. Sturm, \emph{On the geometry of metric measure spaces}. I, Acta Math. 196(2006), 65-131.

\re{}[19]\ B. Thomas, \emph{A natural Finsler-Laplace operator}. arXiv: 1104.4326v2.

\re{}[20]\  B. Wu and Y. Xin, \emph{Comparison theorems in Finsler geometry and their applications}. Math. Ann. 337(2007), 177-196.

\re{}[21]\ G. Wang and C. Xia,  \emph{A sharp lower bound for the first eigenvalue on Finsler manifolds}. arXiv: 1112.4401v1.

\re{}[22]\ J. Zhong and H. Yang, \emph{On the estimates of the first eigenvalue of a compact Riemannian manifold}. Sci. Sinica Ser. A27(1984), 1265-1273.

\vspace*{1.5cm}

Songting Yin\\
1. Department of Mathematics, Tongji University, Shanghai, 200092, China;\\
2. Department of Mathematics and Computer Science, Tongling University, Tongling, 244000 Anhui, China\\
E-mail:yst419@163.com\\

Qun He\\
Department of Mathematics, Tongji University, Shanghai, 200092, China\\
E-mail: hequn@tongji.edu.cn\\

Yibing Shen\\
Department of Mathematics, zhejiang University, Hangzhou, 310028 zhejiang, China\\
E-mail: yibingshen@zju.edu.cn

\newpage
\end{CJK*}
\end{document}